\documentclass[secnum,reqno,11pt]{article}
\usepackage[usenames,dvipsnames]{color}
\usepackage[colorlinks,linkcolor=OliveGreen,urlcolor=MidnightBlue,citecolor=Mahogany]{hyperref}%
\usepackage{hyperref}
\usepackage{mathrsfs}
\usepackage{bookmark}
 \usepackage[all]{hypcap}
\usepackage{amssymb, amsmath, latexsym}
\usepackage{layout}
\usepackage{geometry,graphicx,pict2e,authblk}
\usepackage[numbers,sort&compress]{natbib}
\usepackage{enumerate}
\usepackage{latexsym}
\usepackage{amssymb}
\usepackage{amsmath}

\newtheorem{Lem}{Lemma}[section]

\newtheorem{Def}{Definition}[section]

\newcommand{\E}{{\rm  I~\hspace{-1.15ex}E}}

\newcommand{\cH}{{\cal H}}
\newcommand{\ee}{\varepsilon}

\newcommand{\beq}{\begin{equation}}
\newcommand{\eeq}{\end{equation}}

\newtheorem{teo}{Theorem}[section]

\newtheorem{prop}[teo]{Proposition}
\newtheorem{rmk}[teo]{Remark}

\newcommand{\al}{\alpha}
\newcommand{\Om}{\Omega}
\newcommand{\om}{\omega}
\newcommand{\del}{\delta}

\newcommand{\ffi}{\varphi}

\newcommand{\Ht}{\mathcal{H}_T}

\newcommand{\R}{\mathbb{R}}
\newcommand{\red}{\mathbb{R}^d}

\newcommand{\ubar}[1]{\mkern3mu\underline{\mkern-3mu #1\mkern-3mu}\mkern3mu}

\newenvironment{dem}{\text\bf {\it \textbf{Proof.}}}{$\Box$}
\title{Large deviations for a fractional stochastic heat equation in spatial dimension $\mathbb{R}^{d}$ driven by a spatially correlated noise
\footnotetext{This paper has been accepted for publication in \textit{Stochastics \& Dynamics, Doi: 10.1142/S0219493716500015, 
URL: \url{http://www.worldscientific.com/doi/abs/10.1142/S0219493716500015}}. This reprint differs from the original in
pagination and typographic detail.}}

\date{ \today }

\author{\textsc{El Mellali Tarik\thanks{Supported by the CNRST ``Centre National pour la Recherche Scientifique et Technique", grant $N^\circ$ G03/015, Rabat,
 Morocco.
 Part of this work has been done during a visit, supported by the Action Int\'egr\'ee PHC Volubilis
 MA/10/224, of the Laboratory MAP5 at University Paris-Descartes. We would like to thank this Center for hospitality. } 
 \, and\, Mohamed Mellouk} }
\begin{document}
\maketitle
\begin{abstract}
In this paper we study the Large Deviation Principle (LDP in abbreviation) for  a class of Stochastic Partial Differential Equations 
(SPDEs) in the whole space $\mathbb{R}^d$, with arbitrary dimension $d\geq 1$,
under random influence which is a Gaussian noise, white in time and correlated in space. The
differential operator is a fractional derivative operator. We prove a large deviations principle for our equation,
using a weak convergence approach based on a variational representation of functionals of infinite-dimensional 
Brownian motion. This approach reduces the proof of LDP to establishing basic qualitative properties for controlled
analogues of the original stochastic system.
\end{abstract}
%\newpage
%%
\section{Introduction {and general framework}}
In this paper we consider the following Stochastic Partial Differential
Equation (SPDE in abbreviation) given by
\begin{equation}
\left\{
\begin{array}{l}
\dfrac{\partial u^\varepsilon}{\partial t}(t,x)=\mathcal{D}_{\ubar{\del} }^{\ubar{\al}}
u^\varepsilon(t,x)+b(u^\varepsilon(t,x))+ \sqrt{\varepsilon}\sigma (u^\varepsilon(t,x)){\dot{F}}(t,x), 
\\
u^\ee(0,x)=0,%
\end{array}%
\right.  \label{1}
\end{equation}%
where $\ee>0$,  $(t,x)\in [0,T]\times \mathbb{R}^{d}$, $d\geq 1$, $\ubar{\al} =(\al
_{1},\ldots ,\al _{d})$, $\ubar{\del} =(\del _{1},\ldots ,\del _{d})$ and we will assume, along this paper, that
  $\al _{i}\in ]0,2]\setminus \{1\}\text{ and } |\del _{i}|\leq \min
\{\al _{i},2-\al _{i}\}\text{, }i=1,\ldots ,d$.\medskip\\
$\dot{F}$ is the ``formal" derivative of the Gaussian perturbation and 
$\mathcal{D}_{\ubar{\del } }^{\ubar{\al} }$ denotes a non--local fractional
differential operator on $\red$ defined by
\begin{equation*}
\mathcal{D}_{\ubar{\del} }^{\ubar{\al}}=\sum_{i=1}^{d}D_{\del _{i}}^{\al _{i}},
\end{equation*}%
where $D_{\del _{i}}^{\al _{i}}$ denotes the fractional differential
derivative w.r.t the $i$--th coordinate defined via its Fourier transform $\mathcal{F}$ by
\begin{equation*}
\mathcal{F}\left( D_{\del_i }^{\al_i }\varphi \right) (\xi )=-|\xi
|^{\al_i }\exp \left( -\imath \,\del_i \frac{\pi }{2}\mathrm{sgn}(\xi
)\right) \mathcal{F}\left( \varphi \right) (\xi ),
\end{equation*}%
where $\imath ^{2}+1=0$.\medskip \\ 
\noindent The noise $F(t,x)$ is a martingale measure (in the sense given by Walsh in
\cite{Wa}) to be defined with more details in the sequel. The coefficients $%
b $ and $\sigma :\mathbb{R}\rightarrow \mathbb{R}$ are given functions. We
shall refer to equation \eqref{1} as $Eq_{\ubar{\del},\ee }^{\ubar{\al}} (d,b,\sigma )$.\medskip\\
\noindent The theory of large deviations, which has been investigated for different systems
in recent years, reveals important aspects of asymptotic dynamics. Particular attention has been paid to studying LDP for stochastic
differential equations (SDE) (e.g. \cite{A}, \cite{So}, \cite{C-M},  \cite{Pe}, \cite{F-W2}).\medskip\\
\noindent In this contribution to the theory of Large deviations, first investigated by Freidlin and Wentzell in there original work \cite{F-W} 
for Brownian noise driven SDEs in finite dimension,
we are interested in the study of a stochastic heat equation in spatial dimension $d\geq 1$ driven by 
a spatially correlated noise using the approach in \cite{F-W} (see also \cite{D-Z}).\medskip\\
\noindent In various papers about the LDP for solutions to SPDEs or to stochastic 
evolution equations in a semi-linear framework \cite{bms,Br,B-E-M,B-D,D-Z,D-M,F-W2,MiSo99,Pe}, the strategy used is
similar to the classical one for diffusion processes and based on standard LDP discretization methods which, for a large 
number of problems of significant interest, necessitate exponential probability estimates and
exponential tightness much hardest and most technical.\medskip\\
There exist two distinct methods, in establishing LDP for an SDE or SPDE with multiplicative noise,
the classical approach and the weak convergence method. In
the classical method one should discretize the time horizon and freeze the diffusion
term on each interval and then use the Varadhan's contraction principle. In this
method we should overcome many difficult inequalities for convolution integrals.
In the weak convergence approach, which is the approach employed in this work,
we should obtain some sort of continuity w.r.t. some control variables. We will
clarify further this approach in the current and next section.
Several authors have studied LDP for infinite dimensional SDEs with L\'evy
noise, see \cite{So}, \cite{C-M} for the classical one and \cite{B-C-D,Bo-D} for the weak convergence approach.\medskip\\
In this paper we will prove a tantamount argument of the large deviations
principle, the Laplace principle, and we will study the Uniform Laplace Principle. The reader
should refer to \cite{E-D} for a proof of the aforementioned equivalence.\medskip\\
In order to prove our LDP result, we use a combination of a variational representation for infinite-dimensional
Brownian motion and a transfer principle via Laplace Principle based on compactness and weak convergence proved in \cite{B-D}.
Using this method, an LDP to a reaction-diffusion system has been obtained in 
\cite{B-D-M}. For the case of wave equation in spatial dimension $d=3$, the same method has been applied to derive an LDP result 
in \cite{S-L}, and also in several subsequent papers, for instance see \cite{D-M,S-S}.
The case of heat equation, governed by the same noise, has been considered in \cite{MC-MS} taking the spatial coordinate $x\in[0,1]^d$, $d\geq1$ 
where the authors needs to establish precise estimates of the fundamental solution in order to obtain a Freidlin-Wentzell type inequality.
In contrast, the approach we take in this paper is different from that of \cite{MC-MS}, and it is based on weak
convergence arguments.
We can also refer to original LDP result for the case of a one-dimensional heat equation 
driven by a Brownian sheet in \cite{C-M}.\medskip\\
A summary of the paper is as follows. The first and second subsections will give definitions and preliminary results 
about the fractional operator and the noise both considered in this paper. Next, we present the set of assumptions
required for the existence, uniqueness to equation  $Eq_{\ubar{\del},\ee }^{\ubar{\al}} (d,b,\sigma )$, integrability 
condition for the stochastic integral and also regularity properties of the solution.
The last section is devised into tree subsections.\medskip\\
 The first one reviews basic definitions of LDP, the Laplace principle (LP), the equivalence in between in the 
 setting of Polish spaces and a brief description of the method we follow in this work.
 The statement of the main result along with the set of assumptions to be checked are given in the second subsection.
 The last one will present the blokes of our proofs.
 The value of the constants along this article may change from line to line and somtimes we shall emphasize the
 dependence of these constants upon parameters. 
\subsection{The operator \texorpdfstring{$\mathcal{D}_{\ubar{\del}}^{\ubar{\al}}$}{D}}
In one space dimension, the operator $D_{\del}^{\al}$ is a closed, densely defined operator
on $L^{2}(\mathbb{R})$ and it is the infinitesimal generator of a semi-group
which is in general not symmetric and not a contraction. It is self adjoint
only when $\del=0$ and in this case, it coincides with the fractional
power of the Laplacian.\medskip\\
\noindent According to \cite{Ko,D-I}, $D_{\del}^{\al}$ can be represented for $%
1<\al <2$, by
\begin{equation*}
D_{\del}^{\al}\varphi (x)=\int_{-\infty }^{+\infty }\frac{\varphi
(x+y)-\varphi (x)-y\varphi ^{\prime }(x)}{|y|^{1+\al}}\left( \kappa
_{-}^{\del}\mathbf{1}_{(-\infty ,0)}(y)+\kappa _{+}^{\del}\mathbf{1}%
_{(0,+\infty )}(y)\right) dy,
\end{equation*}%
and for $0<\al<1$, by%
\begin{equation*}
D_{\del}^{\al}\varphi (x)=\int_{-\infty }^{+\infty }\frac{\varphi
(x+y)-\varphi (x)}{|y|^{1+\al}}\left( \kappa _{-}^{\del}\mathbf{1}%
_{(-\infty ,0)}(y)+\kappa _{+}^{\del}\mathbf{1}_{(0,+\infty )}(y)\right)
dy,
\end{equation*}%
where $\kappa _{-}^{\del }$ and $\kappa _{+}^{\del}$ are two
non--negative constants satisfying $\kappa _{-}^{\del}+\kappa
_{+}^{\del}>0$ and $\varphi $ is a smooth function for which the integral
exists, and $\varphi ^{\prime }$ is its derivative. This representation
identifies it as the infinitesimal generator for a non--symmetric $\al$%
--stable L\'{e}vy process.\medskip\\
\noindent Let $G_{\al,\del}(t,x)$ denotes the fundamental solution of the
equation $Eq_{\del,1}^{\al}(1,0,0)$ that is, the unique solution of the
Cauchy problem
\begin{equation*}
\left\{
\begin{array}{l}
\dfrac{\partial u}{\partial t}(t,x)=D_{\del}^{\al}u(t,x), \\\\
u(0,x)=\del_{0}(x),\ \ t>0,\,x\in \mathbb{R},%
\end{array}%
\right.
\end{equation*}%
where $\del _{0}$ is the Dirac distribution. Using Fourier's calculus one get
\begin{equation*}
G_{\al,\del}(t,x)=\frac{1}{2\pi }\int_{-\infty }^{+\infty }\exp
\left( -\imath zx-t|z|^{\al}\exp \left( -\imath \del\frac{\pi }{2}%
\mathrm{sgn}(z)\right) \right) dz.
\end{equation*}
The relevant parameters, $\al$, called the index of \textit{stability}
and $\del$ (related to the asymmetry) improperly referred to as the
\textit{skewness} are real numbers satisfying $\al\in ]0,2]$ and $|\del|\leq \min (\al,\,2-\al).$\newpage
\noindent The function $G_{\al,\del}(t,\cdot )$ has the following properties
(see e.g. \cite{DEB,Dedo}).
\begin{Lem}
\label{lqr} For $\al\in (0,2]\setminus \{1\}$ such that $|\del|\leq
\min (\al,2-\al)$
\begin{description}
\item[$(i)$] The function $G_{\al,\del}(t,\cdot )$ is the density of
a L\'{e}vy $\al$--stable process in time $t$.
\item[$(ii)$] The function $G_{\al,\del}(t,x)$ is not in general
symmetric relatively to $x$.
\item[$(iii)$] Semi-group property: $G_{\al,\del}(t,x)$ satisfies the
Chapman Kolmogorov equation, i.e. for $0<s<t$%
\begin{equation*}
G_{\al,\del}(t+s,x)=\int_{-\infty }^{+\infty }G_{\al  ,\del}
(t,y)G_{\al,\del}(s,x-y)dy.
\end{equation*}
\item[$(iv)$] Scaling property: $G_{\al,\del  }(t,x)=t^{-1/\al}G_{\al,\del}(1,t^{-1/\al}x)$.
\item[$(v)$] There exists a constant $c_{\al}$ such that $0\leq
G_{\al  ,\del }(1,x)\leq c_{\al  }/(1+|x|^{1+\al  })$, for all $x\in
\mathbb{R}$.
\end{description}
\end{Lem}
\noindent Now, for higher dimension $d\geq1$ and any multi index $\ubar{\al } =(\al _{1},\ldots ,\al _{d})$
and $\ubar{\del } =(\del _{1},\ldots ,\del _{d})$,  let $\mathbf{G}_{\ubar{\al } ,\ubar{\del } }(t,x)$
be the Green function of the
deterministic equation $Eq_{\ubar{\del },1 }^{\ubar{\al } }(d,0,0)$. Clearly
\begin{eqnarray*}
\mathbf{G}_{\ubar{\al } ,\ubar{\del } }(t,x) &=&\prod_{i=1}^{d}G_{\al _{i},\del
_{i}}(t,x_{i}) \\
&=&\frac{1}{(2\pi )^{d}}\int_{\mathbb{R}^{d}}\exp \left( -\imath
\left\langle \xi ,x\right\rangle -t\sum_{i=1}^{d}|\xi _{i}|^{\al
_{i}}\exp \left( -\imath \del _{i}\frac{\pi }{2}\mathrm{sgn}(\xi
_{i})\right) \right) d\xi ,
\end{eqnarray*}%
where $\left\langle \cdot ,\cdot \right\rangle $ stands for the inner
product in $\mathbb{R}^{d}$.
\subsection{The driving noise \texorpdfstring{$F$}{F} and a family of SPDEs driven by \texorpdfstring{$F$}{F}}
Let us explicitly describe here the spatially homogeneous
noise (see e.g. \cite{Da}). Precisely, let $\mathcal{S}(\mathbb{R}^{d+1})$ be the space of Schwartz test functions, 
on a complete probability space $(\Omega ,\mathcal{G},P)$, the noise $%
F=\{F(\varphi ),\varphi \in \mathcal{S}(\mathbb{R}^{d+1})\}$ is assumed to
be an $L^{2}(\Omega ,\mathcal{G},P)$--valued Gaussian process with mean zero
and covariance functional given by
\begin{equation*}
J(\varphi ,\psi ):=\E(F(\varphi)F(\psi) )=\int_{\mathbb{R}_{+}}\int_{\mathbb{R}^{d}}
\left( \varphi (s,\star )\ast \widetilde{\psi }(s,\star)\right)(x)\Gamma
(dx) ds,\, \,\varphi ,\psi \in \mathcal{S}(\mathbb{R}^{d+1}),
\end{equation*}%
where $\widetilde{\psi }(s,x)=\psi (s,-x)$ and $\Gamma $ is a non--negative
and non--negative definite tempered measure, therefore symmetric.
The symbols $*$ denotes the convolution product and $\star$ stands for the spatial variable.\medskip\\
\noindent Let $\mu $
denotes the spectral measure of $\Gamma $ (usually called the spectral measure of the
noise $F$), which is also a trivial tempered measure
(see \cite[Chap. VII, Th\'eor\`eme XVII]{Sch}), that is $\mu =\mathcal{F}^{-1}(\Gamma )$ and this gives
\begin{equation}
J(\varphi ,\psi )=\int_{\mathbb{R}_{+}}ds\int_{\mathbb{R}^{d}}\mu (d\xi )%
\mathcal{F}\varphi (s,\cdot )(\xi )\overline{\mathcal{F}\psi (s,\cdot )}(\xi
),  \label{cov}
\end{equation}%
where $\overline{z}$ is the complex conjugate of $z$.\medskip\\
\noindent A typical example of space correlation is given by $\Gamma(x)=f(x)dx$, where $f$ is a nonnegative
function which is assumed to be integrable around the origin. In this case, the
covariance functional $J$ reads
$$J(\varphi ,\psi )=\int_{\mathbb{R}_{+}}\int_{\mathbb{R}^{d}} \int_{\mathbb{R}^{d}} \varphi(s,x )f(x-y)\psi(s,y) dy\, dx\, ds\,.$$
The space-time white noise would correspond to the case where $f$ is the Dirac delta at the origin.\\
\noindent Following the same approach as in \cite{Da}, the Gaussian process $F$
can be extended to a worthy martingale measure $\{F(t,A):=F([0,t]\times
A)\,;\,t\geq 0,\,A\in \mathcal{B}_{b}(\mathbb{R}^{d})\}$ which shall acts as
integrator, in the sense of Walsh \cite{Wa}, where $\mathcal{B}_{b}(\mathbb{R%
}^{d})$ denotes the bounded Borel subsets of $\mathbb{R}^{d}$. Let $\mathcal{%
G}_{t}$ be the completion of the $\sigma $--field generated by the random
variables $\{F(s,A);\;0\leq s\leq t,\;A\in \mathcal{B}_{b}(\mathbb{R}^{d})\}$%
. The properties of $F$ ensure that the process $F=\{F(t,A);\;t\geq 0,\;A\in
\mathcal{B}_{b}(\mathbb{R}^{d})\}$, is a martingale with respect to the
filtration $\{\mathcal{G}_{t};\,t\geq 0\}$.\medskip\\
\noindent Then, one can give a rigorous meaning to the solution of equation $Eq_{\ubar{\del },\ee
}^{\ubar{\al } }(d,b,\sigma )$, by means of a jointly measurable and $\mathcal{G}%
_{t}$--adapted process $\{u(t,x);(t,x)\in \mathbb{R}_{+}\times \mathbb{R}%
^{d}\}$ satisfying, for each $t\geq 0$ and a.s. for almost all $x\in \mathbb{R}^{d}$ the following evolution equation:
\begin{eqnarray}
u^\ee(t,x) &=&\sqrt{\ee}\int_{0}^{t}\int_{\mathbb{R}^{d}}\mathbf{G}_{\ubar{\al } ,\ubar{\del }
}(t-s,x-y)\sigma (u^\ee(s,y))F(ds,dy) \notag \\
&&+\,\int_{0}^{t}ds\int_{\mathbb{R}^{d}}\mathbf{G}_{\ubar{\al } ,\ubar{\del }
}(t-s,x-y)b(u^\ee(s,y))dy.  \label{3}
\end{eqnarray}%
\noindent Note that the last  stochastic integral on the right-hand side of \eqref{3} can be understood in the sense of Walsh \cite{Wa} 
or using the further extension of Dalang \cite{Da}.\medskip \\
\noindent In order to prove our main result, we are going to give another equivalent 
approach to the solution of  $Eq_{\ubar{\del},\ee }^{\ubar{\al}} (d,b,\sigma )$, see \cite{Dg-Qr}.
To start with, let us denote by $\cH$ the Hilbert space obtained by the completion of
$\mathcal{S}(\R^d)$ with the inner product
$$\langle\ffi,\psi\rangle_{\cH}=\int_{\R^d}\Gamma(dx)\,
(\ffi\ast\tilde{\psi})(x)=\int_{\mathbb{R}^{d}}\mu (d\xi )%
\mathcal{F}\varphi(\xi )\overline{\mathcal{F}\psi}(\xi
);\,\ffi,\,\psi\in\mathcal{S}(\R^d).$$
\noindent By the Walsh theory of martingale measures \cite{Wa}, for $t\geq0$ and  $h\in\cH$
 the stochastic integral $$B_t(h) = \int_0^t \int_{\R^d} h(y) F(ds,dy),$$ is
well defined and the process $\{B_t(h);\, {t\geq0},  h\in\cH\}$  is a cylindrical Wiener process on $\cH$, that is: (a) For every 
$h \in \cH$ with  $\|h\|_\cH=1,$  $\{B_t(h)\}_{t\geq 0}$ is a standard Wiener process;
and (b) For every $t\geq 0, a,b \in \R$ and $f,g \in \cH,$ $B_t(af+bg)=aB_t(f)+bB_t(g)$ almost surely.\medskip\\
\noindent  Let  $(e_k)_{k\geq1}$  be
a complete orthonormal system (\textit{CONS}) of the Hilbert
space $\cH$, then $$\{B_t^k:=\int_0^t\int_{\R^d}e_k(y)F(ds,dy);\,k\geq1\}$$  defines a sequence of
independent standard Wiener processes and we have the following representation
\beq
B_t=\sum_{k\geq1}B_t^k e_k.
\eeq
\noindent  Let $(\mathscr{F}_t)_{t\in[0, T]}$ be the $\sigma$-field generated by the random variables 
$\{B_s^k;\, s\in[0,t],\,k\geq1\}$. We define the predictable $\sigma$-field as the $\sigma$-field in $\Omega \times [0,T]$ 
generated by the sets $\{(s,t]\times A;\,A\in\mathscr{F}_s,\,0\leq s<t\leq T\}$. 
In the following we can define the stochastic integral with respect to cylindrical Wiener process
$(B_t(h))_{t\geq0}$ (see e.g. \cite[Chapter 4]{Dap} or \cite{Dg-Qr}) of any predictable square-integrable process with values in $\cH$ as follows
$$\int_0^t \int_{\R^d}g\cdot dB:=\sum_{k\geq1}\int_0^t\,\langle g(s), e_k\rangle_{\cH}\,dB_s^k.$$
Note that the above series converge in $L^2(\Omega,\mathscr{F}, P)$ and the sum does not depend on the selected \textit{CONS}.
Moreover,  each summand, in the above series, is a classical It\^o integral with respect to a standard Brownian motion, 
and the resulting stochastic integral is a real-valued random variable.\medskip\\
\noindent In the sequel, we shall consider the $\it{ mild }$ solution to equation $Eq_{\ubar{\del},\ee }^{\ubar{\al}} (d,b,\sigma )$
given by
\begin{eqnarray}
\label{4}
 u^{\ee}(t,x) &=& \sqrt{\ee}\sum_{k\geq1}\int_0^t\,\langle \mathbf{G}_{\ubar{\al } ,\ubar{\del } }(t-s,x-\cdot)\sigma(u^{\ee}(s,\cdot)),e_k\rangle_{\cH}\,dB_s^k \\
&& \,\, \quad \,\, \quad  \,\, \quad +  \int_0^t\left[\mathbf{G}_{\ubar{\al } ,\ubar{\del } }(t-s)\ast b(u^{\ee}(s,\cdot))\right](x)ds,\nonumber
\end{eqnarray}
$t\in[0,T]$, $x\in\red$ and $``\ast "$ stands for the convolution operator.
\subsection{Existence, uniqueness and H\"older regularity to equation \texorpdfstring{$Eq_{\ubar{\del},\ee}^{\ubar{\al}}(d,b,\sigma)$}{Eq. (1.1)}}
The purpose of this section is to give sufficient conditions for the existence
and uniqueness to our equation and also H\"older regularity of the solution which we will use in the sequel.
\noindent Now, for a given multi index $\ubar{\al }=(\al_{1},\dots,\al_{d})$ such
that $\al_{i}\in]0,2]\setminus\{1\}$, $i=1,\dots,d$ and any $\xi\in\red$, let
\begin{equation*}
S_{\ubar{\al }}(\xi)=\sum_{i=1}^{d}|\xi_i|^{\al_i}.
\end{equation*} 
Assume the following assumptions on the functions $\sigma$, $b$ and the
measure $\mu:$
\begin{description}
 \item[(C)]: The functions $\sigma$ and $b$ are Lipschitz.
\item[$\mathbf{(H_{\eta}^{\ubar{\al }})}$]: Let $\ubar{\al}$ as defined above and $\eta\in(0,1]$
 $$\int_{\red}\frac{\mu(d\xi)}{(1+S_{\ubar{\al }}(\xi))^{\eta}}<\infty.$$
\end{description}
\noindent The last assumption stands for an integrability condition w.r.t. the spectral measure $\mu$. Indeed, the following stochastic integral
\begin{equation*}
\int_{0}^{T}\int_{\red}\mathbf{G}%
_{\ubar{\al } ,\ubar{\del } }(T-s,x-y)F(ds,dy),
\end{equation*}
is well defined if and only if 
\begin{equation*}
\int_{0}^{T}ds\int_{\mathbb{R}^{d}}\mu (d\xi )|{\mathcal{F}\mathbf{G}%
_{\ubar{\al } ,\ubar{\del } }}(s,\cdot )(\xi )|^{2}<+\infty .
\end{equation*}
More precisely, appealing to \cite[Lemma 1.2]{B-E-M} there exist two positive constants $c_{1}$ and $c_{2}$ such
that
\begin{equation}
\label{eqfourie}
c_{1}\int_{\mathbb{R}^{d}}\frac{\mu (d\xi )}{1+S_{\ubar{\al } }(\xi )}\leq
\int_{0}^{T}\int_{\mathbb{R}^{d}}\mu (d\xi )\left\vert \mathcal{F}\mathbf{G}%
_{\ubar{\al } ,\ubar{\del } }(s,\cdot )(\xi )\right\vert ^{2}ds\leq c_{2}\int_{\mathbb{R%
}^{d}}\frac{\mu (d\xi )}{1+S_{\ubar{\al } }(\xi )}.
\end{equation}  
\begin{rmk}
The upper and lower bounds in \eqref{eqfourie} do not depend on the
parameter $\ubar{\del } $. When $\al _{i}=2$ for all $i=1,\ldots ,d$ then $%
S_{2}(\xi )=\sum_{i=1}^d|\xi_i|^2=:\left\vert \xi \right\vert^{2}$ and the bounds in 
\eqref{eqfourie} are the same ones which appeared in \cite{Da} (see also \cite{Sa}), that is
\begin{equation*}
\int_{\mathbb{R}^{d}}\frac{\mu (d\xi )}{1+\left\vert \xi \right\vert^{2}%
}<+\infty.
\end{equation*}
\end{rmk}
\noindent
In the paper \cite{B-D-M}, Theorem 2.1 gives existence and uniqueness of random field solutions to equation \eqref{4} 
under conditions $\mathbf{(C)}$ and $\mathbf{(H_1^{\ubar{\al}})}$ for $\eta=1$.
 In fact, this result extend those of \cite[Theorem 2.1]{AM} \cite[Theorem 1]{Dedo} to the $d-$dimensional case and 
\cite[Theorem 13]{Da} to the fractional setting.\medskip\\
\noindent
Moreover, \cite[Theorem 3.1]{B-E-M} gives the regularity properties of the solution to equation \eqref{4} in time and space
 improving the results in \cite{Hc,Pp} corresponding to the case $\al_i=2$, $\del_i=0$ for $i=1,\dots,d$.
 More precisely, the trajectories of the solution are $\beta=(\beta_1,\beta_2)-$H\"older continuous in $(t,x)\in[0,T]\times K$ for every 
$\beta_1\in(0,\al_0\frac{1-\eta}{2})$, $\beta_2\in(0,1-\eta)$ and every $K$ compact subset of $\red$,
where $\al_0=\min_{1\leq i\leq d}\{\al_i\}$.\medskip\\
\noindent Consequently, the random field solution $\{u(t,x);\,(t,x)\in[0,T]\times K\}$ to equation \eqref{4} lives in the H\"older space
 $C^{\beta}([0,T]\times K;\red)$ equipped with the norm defined by
\begin{equation}\label{N}
\left\Vert f\right\Vert _{\beta,K}=\sup_{(t,x)\in [0,T]\times
K}\left\vert f(t,x)\right\vert +\sup_{s\neq t\in [0,T]}\sup_{x\neq
y\in K}\dfrac{\left\vert f(t,x)-f(s,y)\right\vert }{\left\vert
t-s\right\vert ^{\beta_{1}}+\left\Vert x-y\right\Vert ^{\beta_{2}}}.
\end{equation}
 \noindent Let us define the space $\Ht\stackrel{.}{=}L^2([0,T];\cH)$. For any $h\in\Ht$, consider the deterministic evolution equation
\beq\label{6}
\begin{split}
Z^{h}(t,x)=\int_0^t\,\langle\mathbf{G}_{\ubar{\al},\ubar{\del}}(t-s,x-\cdot)\sigma(Z^{h}(s,\cdot)),
h(s,\cdot)\rangle_{\cH}\,ds\\
+\int_0^t\left[\mathbf{G}_{\ubar{\al},\ubar{\del}}(t-s)\ast b(Z^{h}(s,\cdot))\right](x)ds,
\end{split}
\eeq
where the first term on the right-hand side of the above equation can be written as
\begin{equation*}
 \sum_{k\geq1}\int_0^t\,\langle\mathbf{G}_{\ubar{\al},\ubar{\del}}(t-s,x-\cdot)\sigma(Z^{h}(s,\cdot)),e_k\rangle_{\cH}h_k(s)\,ds,
\end{equation*}
with $h_k(t)=\langle h(t),e_k\rangle_{\cH}$, $t\in[0,T]$, $k\geq1$.\\
For the existence and H\"olderian regularities for the solution of equation \eqref{6} see Remark \ref{RM}.
\section{General framework and Large deviation principle result}
\subsection{Large deviation principle and Laplace principle}
Let  $\{X^\ee;\, \ee>0\}$ be a family of random variables defined on a probability space 
$(\Omega, \mathscr{F}, P)$ and taking values in a Polish space (i.e. separable complete metric space) ${\cal E}$.
We denote by $\E$ the expectation with respect to $P$ . The LDP for the family $\{X^\ee;\, \ee>0\}$ 
is concerned with events $A$ for which probabilities $P(X^\ee \in A)$ converges to zero exponentially fast as $\ee \rightarrow 0$.
The exponential decay rate of such probabilities are typically expressed in terms of a ``rate function" $I$ mapping ${\cal{E}}$ into $[0, \infty]$. 
 \begin{Def}
 The family of random variables $\{X^\ee;\, \ee>0\}$ is said to satisfy the LDP with the good rate function (or action functional)
 $I: {\cal{E}} \rightarrow  [0, \infty]$,  on ${\cal E}$, if 
 \begin{enumerate}
 \item For each $M<\infty$ the level set $\{x\in  {\cal{E}};\,  I(x) \leq M\}$ is a compact subset of  $\cal{E}$. 
 \item \it{Large deviation upper bound:} for any closed subset $F $ of $\cal{E}$
$$\limsup_{\ee\rightarrow0^+}\ee\log P(X^{\ee}\in F)\leq -I(F).$$
\item \it{Large deviation lower bound:} for any open subset $O$ of  $\cal{E}$
$$\liminf_{\ee\rightarrow0^+}\ee\log P(X^{\ee}\in O)\geq -I(O).$$
\end{enumerate}
Where, for $A \subset\cal{E}$, we define $I(A)=\inf_{x \in A}I(x)$.
 \end{Def}
\noindent Varadhan's and Bryc's results, \cite{V} and \cite{Br}, announced an equivalence between
LDP and Laplace principle (LP), which notices the expectations of exponential
functions.
\begin{Def}(Laplace principle)
 The family of random variables $\{X^\ee;\, \ee>0\}$ defined on the Polish space ${\cal E}$, is said to satisfy the Laplace principle with the
rate function $I$ 
 if for any bounded continuous function $h:{\cal E} \rightarrow \R$,
{$$\lim_{\ee\rightarrow 0}\ee\log \E\left( \exp\left[-\frac{1}{\ee}h(X^\ee)\right]\right)= -\inf_{f\in {\cal E}}\{h(f)+I(f)\}.$$}
  \end{Def}
Another display of variational representation in evaluating the exponential
integrals is in the following proposition which is a cornerstone of weak convergence
method. For a comprehensive introduction to the applications of weak
convergence method to the theory of large deviations we refer the reader to the paper \cite{E-D}.
\begin{prop}
Let $({\cal V}, {\cal A})$ be a measurable space and $f$ be a bounded
 measurable function mapping ${\cal  V}$  into the real numbers $\R$. For a given probability
measure $ \theta$ on  ${\cal  V}$, we have the following representation
$$-\log \int_{\cal  V} e^{-f} d\theta= \inf_{\gamma \in {\cal P} ({\cal V})}\{R(\gamma ||\theta)+\int_{\cal  V}f d \gamma\},$$
where
$R(\gamma ||\theta):=\int_{\cal  V} \log (\frac{d \gamma}{d \theta}) d\gamma$ and $ {\cal P} ({\cal V})$ 
denotes the set of  probability measures on ${\cal V}$.
\end{prop}
By using the above Proposition, the following variational representation has been obtained in \cite{B-D,B-D-M}
for exponential integrals w.r.t. the Wiener process.
\subsubsection*{Variational representation}
Let $B=\{B_k(t);\,t\in[0,T],\,k\geq1\}$ be a sequence of independent real
standard Brownian motions, and notice that $B$ is a
$C( [0,T]; \mathbb{R}^\infty)$-valued 
random variable. Consider the Hilbert space 
$l_2\dot{=}\{x\equiv(x_1,x_2,...)\,;\,x_i\in\R\,\text{and}\,\sum x_i^2<\infty\}$,
and let $\mathcal{P}_2(l_2)$ be the family of all $l_2$-valued
predictable processes for which $\int_0^T\left\|\phi(s)\right\|^2_{l_2}ds<\infty$ a.s., 
where $\|\cdot\|_{l_2}$ denotes the norm in the Hilbert space $l_2$.
That is, $u\in\mathcal{P}_2(l_2)$ can be written as $u=\{u_i\}^{\infty}_{i=1}$, 
$u_i\in\mathcal{P}_2(\R)$ and $\sum_{i=1}^{\infty}\int_0^T|u_i(s)|^2ds<\infty$ a.s.
\begin{teo}(\cite[Theorem 2]{B-D-M}).
Let $g$ be a bounded, Borel measurable function mapping $C([0,T];\,\mathbb{R}^{\infty})$ into $\R$. Then
\begin{equation*}
-\log\E\left[\exp(-g(B))\right]=\inf_{u\in{\cal{P}}_2(l_2)}\E\left[\frac{1}{2}\int^T_0\|u(s)\|_{l_2}^2ds+g\left(B+\int^{.}_0u(s)ds\right)\right].
 \end{equation*}
\end{teo}
According to \cite{Bo-D}, whenever the functionals of interest are expressed as measurable functionals 
of a Wiener process, the above stated representation can be used to derive various asymptotic results 
of large deviations nature.\medskip\\
\noindent This representation together with Laplace's Principle present a different method in obtaining LDP for large class
of stochastic equation driven by a Gaussian type noise, by using stochastic control and weak convergence approach
for a given family ${\cal G}^{\ee}(B(\cdot))$, where ${\cal G}^{\ee}$ is an appropriate family of measurable maps from the Wiener space
to some Polish space and $B(\cdot)$ stands for a Hilbert space valued Wiener process (see \cite{E-D}).
\subsection{The main result}
\noindent 
The aim of this work is to apply the weak convergence approach to establish
LDP  for the family 
$\{u^{\ee};\,\ee\in(0,1]\}$ given by \eqref{4}, in the space of $\beta$--H\"older continuous $C^{\beta}([0,T]\times K;\red)$,
with the rate function $I$ defined below by \eqref{I}.\medskip\\
\noindent As mentioned in the Introduction section, to put our problem in its obvious setting, we need a Polish space carrying the probability laws 
of the family $\{u^{\ee}(t,x);\,\ee\in(0,1],\,(t,x)\in[0,T]\times\red\}$. Since $C^{\beta}([0,T]\times K;\red)$ 
is not separable, we are brought to consider the space $C^{\beta',0}([0,T]\times K;\red)$ of H\"older continuous functions $f$
with degree $\beta'<\beta$ such that 
$$\lim_{\del\rightarrow0^+}\left(\sup_{|t-s|+|x-y|<\del}\frac{|f(t,x)-f(s,y)|}{|t-s|^{\beta'_1}+\|x-y\|^{\beta'_2}}\right)=0,$$
and $C^{\beta',0}([0,T]\times K;\red)$ is a Polish space containing $C^{\beta}([0,T]\times K;\red)$.
\medskip\\%%
From now on, let $\mathcal{E}^{\beta}:=C^{\beta,0}([0,T]\times K;\red)$ be the space of
$(\beta_1,\beta_2)$-H\"{o}lder continuous functions in time and space
respectively which we equip with the 
H\"older norm of degree $\beta$ defined by \eqref{N}, where 
$\beta=(\beta_1,\beta_2)$ and satisfying $0<\beta_1<\frac{\al_0(1-\eta)}{2}$, $0<\beta_2<1-\eta$ and
$\al_0=\min_{1\leq i\leq d}\{\al_i\}$.\medskip\\
\noindent We introduce the map $\mathcal{G}^0$ that will be used to define the rate function in our setting, that is
\begin{align}\label{g0}
\begin{aligned}
 \mathcal{G}^0\,:\, &\Ht\longrightarrow {\cal{E}}^{\beta}\\
                 & h \longmapsto \mathcal{G}^0(h)=Z^h,
\end{aligned}
\end{align}
where $Z^h$ is the strong solution of the integral equation defined by \eqref{6}.
\newpage%%
\noindent Our aim is to prove the following.
\begin{teo}\label{21}
 Assume $\mathbf{(C)}$ and $\mathbf{(H_{\eta}^{\ubar{\al }})}$ for $\eta\in(0,1]$ and $\ubar{\al }=(\al_1,\dots,\al_d)$, 
 $\al_i\in]0,2]\setminus\{1\}$ for $i=1,\dots,d$. Let 
$\{u^{\ee}(t,x);\,(t,x)\in[0,T]\times \red\}$ be the solution of 
equation \eqref{4}. Then, the law of the solution $\{u^{\ee};\,\ee\in(0,1]\}$ satisfies on
${\cal{E}}^{\beta}$, a large deviation principle 
with rate function
\begin{equation}\label{I}
 I(f)=\inf_{\{h\in\Ht:\,\mathcal{G}^0(h)=f\}}\{\frac{1}{2}\|h\|^2_{\Ht}\},
\end{equation}
where  $\mathcal{G}^0(\cdot)$ is defined by \eqref{g0}.
\end{teo}
In order to prove Theorem \ref{21}, we adopt the weak convergence approach \cite{E-D}.
According to \cite{B-D}, the crucial step in the proof is a variational 
representation for some functionals of an infinite dimensional Brownian 
motion combined with a transfer principle via Laplace Principle based on compactness and weak convergence.
The authors also states that this method can be viewed as an extended contraction 
principle which allows to derive a Wentzell-Freidlin type large deviation results.\medskip\\
\noindent Accordingly, and based on this approach, we consider a set of assumptions that will be used to ensure the validity
of Theorem \ref{21}.
Now we can formulate the following sufficient conditions of the
Laplace principle (equivalently, Large deviation principle) in \cite{B-D} for $u^\ee$ as $\ee \rightarrow 0$.
\subsubsection*{Weak regularity}
Denote by $\mathcal{A}_{\cH}$ the set of predictable processes which belong to
$L^2(\Om\times[0,T];\cH)$. For any $N>0$, define 
$$\Ht^N\stackrel{.}{=}\{h\in\Ht;\,\|h\|_{\Ht}\leq N\},$$
$$\mathcal{A}_{\cH}^N\stackrel{.}{=}\{u(\om)\in\mathcal{A}_{\cH};\,u\in\Ht^N\, a.s.\},$$ 
and we consider that $\Ht^N$ is endowed with the weak topology of $\Ht$.\medskip\\
The sets $\mathcal{A}_{\cH}$ and $\Ht^N$ defined above will play a central role in the proofs of the Laplace principle. Indeed, 
these sets are essential in proving tightness for sequences of Hilbert space valued processes by applying Theorem 2.4 and 
Theorem 2.5 in \cite{B-D}, which we will consider later.\\
\noindent For $v\in\mathcal{A}_{\cH}^N$, $\ee\in(0,1]$ define the controlled equation
$u^{\ee,v}$ by
\begin{eqnarray}
\label{7}
 u^{\ee,v}(t,x)&=&\sqrt{\ee}\sum_{k\geq1}\int_0^t\,\langle\mathbf{G}_{\ubar{\al},\ubar{\del}}(t-s,
x-\cdot)\sigma(u^{\ee,v}(s,\cdot)),e_k\rangle_{\cH}\,dB_s^k\nonumber\\
&&\quad \quad +\int_0^t\,\langle\mathbf{G}_{\ubar{\al},\ubar{\del}}(t-s,x-\cdot)\sigma(u^{\ee,v}(s,\cdot)),v(s,\cdot)\rangle_
{\cH}\,ds \\
&&\quad \quad +\int_0^t\mathbf{G}_{\ubar{\al},\ubar{\del}}(t-s)\ast b(u^{\ee,v}(s,\cdot))](x)ds.\nonumber
\end{eqnarray}
\newpage
\noindent Then, let's consider the following two conditions which correspond to the weak convergence approach framework in our setting
\begin{description}
 \item[a)] The set $\{Z^h;\, h\in\Ht^N\}$ is a compact set of
${\cal{E}}^{\beta}$, $Z^h$ being the solution of equation
\eqref{6}.
\item[b)] For any family $\{v^{\ee};\,\ee>0\}\subset\mathcal{A}_{\cH}^N$ which
converges in distribution as $\ee\rightarrow0$ to 
$v\in\mathcal{A}_{\cH}^N$, as $\Ht^N$-valued random variables, we
have $$\lim_{\ee\rightarrow0}u^{\ee,v^{\ee}}=Z^v$$
in distribution, as ${\cal{E}}^{\beta}$-valued random variables, where $Z^v$ denotes the solution of \eqref{6} corresponding to the  $\Ht^N$-valued random
variable $v$ (instead of a deterministic function $h$).
\end{description}
\begin{rmk}
We should at this point give a commentary about the two conditions considered above. Condition $\mathbf{a)}$ says that the level sets 
of the rate function are compact, and condition $\mathbf{b)}$ is a crucial assumption in the application of the weak convergence approach 
and is a statement of weak convergence of the family of random variables $\{u^{\ee,v^{\ee}};\,\ee>0\}$ as $\ee$ goes to 0. 
\end{rmk}
\noindent Let  $$\mathcal{G}^\ee: C([0,T];\mathbb{R}^\infty)\longrightarrow {\cal E}^{\beta}, \, \ee>0$$ be a family
of measurable maps such that $\mathcal{G}^{\ee}(B(\cdot)):=u^\ee$
 (where $u^{\ee}$ stands for the solution to equation \eqref{4}).\medskip\\
\noindent
Then, by applying \cite[Theorem 6]{B-D-M} to the above defined functional $\mathcal{G}^{\ee}$ and 
 $\mathcal{G}^0$ given by \eqref{g0}, a verification of conditions \textbf{a)} and \textbf{b)} implies 
the validity of Theorem \ref{21}.
\subsection{Proof of Theorem \ref{21}}
Both conditions $\mathbf{a)}$ and $\mathbf{b)}$ will follow
from a single continuity result. Condition $\mathbf{a)}$ will follow by proving 
the continuity of the mapping $h:\Ht^N\rightarrow Z^h\in {\cal E}^{\beta}$
 with respect to the weak topology.\medskip\\
 It will consist on proving that, if for 
 $h$, $(h_n)_{n\geq1}\subset\Ht^N$ such that for any $g\in\Ht$,
$$\lim_{n\rightarrow\infty}\langle h_n-h,g\rangle_{\Ht}=0,$$
then,
\beq\label{8}
\lim_{n\rightarrow\infty} \|Z^{h_n}-Z^h\|_{\beta,K}=0
\eeq 
For the condition $\mathbf{b)}$, we will use the Skorohod representation theorem to reformulate it.
That is, there exist a probability space $(\overline{\Om},\overline{\mathscr{F}},\overline{P})$, a sequence of independent Brownian motions
$\overline{B}=(\overline{B}_k)_{k\geq1}$ along with the 
corresponding filtration $(\overline{\mathscr{F}}_t)_{t\in[0,T]}$ where
$\overline{\mathscr{F}}_t=\sigma\{\overline{B}_k(s);\,0\leq s\leq t,\,k\geq1\}$ and also a family of $\overline{\mathscr{F}}_t$--predictable processes
$(\overline{v}^{\ee},\,\ee>0)$, $\overline{v}$ belonging to 
$L^2(\overline{\Om}\times[0,T];\mathcal{H})$ taking values on $\Ht^N$
$\overline{P}$ a.s., such that the joint law of
$(v^{\ee},v,B)_{P}$ 
coincides with that of
$(\overline{v}^{\ee},\overline{v},\overline{B})_{\overline{P}}$ and
satisfying
$$\lim_{\ee\rightarrow0}\langle\overline{v}^{\ee}-\overline{v},g\rangle_{\Ht}
=0,\quad\overline{P}\,a.s.,\quad g\in\Ht,$$
as $\Ht^N$-valued random variables.\medskip\\
Let $\overline{u}^{\ee,\overline{v}^{\ee}}(t,x)$ be the solution to a similar
equation as \eqref{7} obtained 
by changing $v$ into $\overline{v}^{\ee}$ and $B_k$ by $\overline{B}_k$. Thus, verifying condition $\mathbf{b)}$
 will consist on proving that for any $q\in[1,\infty[$ we have
\beq\label{9}
\lim_{\ee\rightarrow0}\overline{\E}\left[\left\|\overline{u}^{\ee,
\overline{v}^{\ee}}-Z^{\overline{v}}\right\|^q_{\beta,K}\right]=0,
\eeq
$\overline{\E}$ being the expectation operator on
$(\overline{\Om},\overline{\mathscr{F}},\overline{P})$.\medskip \\
In the other hand, notice that taking $\ee=0$ and substitute $\overline{v}$ for $h\in\mathcal{A}_{\cH}^N$ in \eqref{7}
we obtain the deterministic evolution equation \eqref{6} satisfied by $Z^h$. Accordingly, the convergence \eqref{8} will 
follow once \eqref{9} is proved.
\noindent According to Lemma A1 in \cite{bms}, the proof of \eqref{9} can be carried into two
steps :
\begin{description}
 \item[1- Estimation of the increments] 
\beq\label{10}
\begin{split}
\sup_{\ee\leq1}\overline{\E}\left(\left|\left[\overline{u}^{\ee,\overline{v}^{
\ee}}(t,x)-Z^{\overline{v}}(t,x)\right]-
\left[\overline{u}^{\ee,\overline{v}^{\ee}}(s,y)-Z^{\overline{v}}(s,y)\right]
\right|^q\right)\\
 \leq C\left[|t-s|^{\beta_1}+\|x-y\|^{\beta_2}\right]^q.
\end{split}
\eeq
\item[2- Point-wise convergence]
 \beq\label{11}
\lim_{\ee\rightarrow0}\overline{\E}\left(\left|\overline{u}^{\ee,\overline{v}^{
\ee}}(t,x)-Z^{\overline{v}}(t,x)\right|^q\right)=0,
\eeq
where $q\in[1,\infty[$, $(t,x),\,(s,y)\in[0,T]\times K$.
\end{description}
\noindent First, we show the following proposition which stands for a statement of existence and uniqueness
of the stochastic controlled equation given by \eqref{7}
\begin{prop}\label{12}
 Assuming $\mathbf{(C)}$ and $\mathbf{(H_{\eta}^{\ubar{\al }})}$, for $\eta\in(0,1]$ and $\ubar{\al }=(\al
_{1},\ldots ,\al _{d})$ satisfying $\al_{i}\in]0,2]\setminus\{1\}$, $i=1,\dots,d$. Then, there exists a unique random field
solution to equation \eqref{7}, $\{\overline{u}^{\ee,\overline{v}^{\ee}}(t,x)
;\,(t,x)\in[0,T]\times\red\}$, which satisfies
\begin{equation}\label{estm}
 \sup_{\ee\leq1}\sup_{v\in\mathcal{A}_{\cH}^N}\sup_{(t,x)\in[0,T]\times\red}
\overline{\E}\left[\left|\overline{u}^{\ee,\overline{v}^{\ee}}(t,x)\right|^q\right]<\infty.
\end{equation}
\end{prop}
\begin{dem}
From now on, we drop the bars in the notation for the sake of simplicity. We only sketch the steps of the proof following
 those of \cite[Theorem 2.1]{B-E-M}, and is based on the Picard iteration scheme
\newline
\beq
\begin{split}\label{P0}
 u^{\ee,v^{\ee}}_{(0)}(t,x)=0& \\
 u^{\ee,v^{\ee}}_{(n+1)}(t,x)=&\sqrt{\ee}\sum_{k\geq1}\int_0^t\,\langle \mathbf{G}_{\ubar{\al},\ubar{\del}}(t-s,
x-\cdot)\sigma(u^{\ee,v^{\ee}}_{(n)}(s,\cdot)),e_k\rangle_{\cH}\,dB_s^k\\
+&\int_0^t\,\langle\mathbf{G}_{\ubar{\al},\ubar{\del}}(t-s,x-\cdot)\sigma(u^{\ee,v^{\ee}}_{(n)}(s,\cdot)),v^{\ee}(s,\cdot)\rangle_
{\cH}\,ds\\
+&\int_0^tG_{\ubar{\al},\ubar{\del}}(t-s)\ast b(u^{\ee,v^{\ee}}_{(n)}(s,\cdot))](x)ds.
\end{split}
\eeq
The first step is to check that the process $\{u^{\ee,v^{\ee}}_{(n)}(t,x);\,(t,x)\in[0,T]\times\red\}$ 
is well-defined and, for $q\geq1$ 
\beq\label{P1}
\sup_{\ee\leq1}\sup_{v^{\ee}\in\mathcal{A}_{\cH}^N}\sup_{(t,x)\in[0,T]\times\red}
\E\left[\left|u^{\ee,v^{\ee}}_{(n)}(t,x)\right|^q\right]<\infty.
\eeq
\noindent Then
\beq\label{P2}
\sup_{n\geq0}\sup_{\ee\leq1}\sup_{v^{\ee}\in\mathcal{A}_{\cH}^N}\sup_{(t,x)\in[0,T]\times\red}
\E\left[\left|u^{\ee,v^{\ee}}_{(n)}(t,x)\right|^q\right]<\infty,
\eeq
that is the bound \eqref{P1} holds uniformly with respect to $n$.\medskip\\
\noindent Secondly, for $n\geq0$ let
\begin{equation*}
M_n(t):=\sup_{(s,x)\in[0,t]\times\red}
\E\left[\left|u^{\ee,v^{\ee}}_{(n+1)}(s,x)-u^{\ee,v^{\ee}}_{(n)}(s,x)\right|^q\right],
\end{equation*}
then, we prove that
\beq\label{P-3}
M_{n+1}(t)\leq C_q\int_0^t M_n(s)\left(1+{\cal J}(t-s)\right)ds .
\eeq
where 
\begin{equation}\label{J}
{\cal J}(t-s)=\int_{\red}\mu(d\xi)\left|\mathcal{F}\mathbf{G}_{\ubar{\al},\ubar{\del}}(t-s)(\xi)\right|^2. 
\end{equation}
Consequently, we can affirm that the sequence $\{u^{\ee,v^{\ee}}_{(n)}(t,x);\,n\geq0\}$ 
converge in $L^q(\Omega)$, uniformly in $(t,x)$, to a limit $u^{\ee,v^{\ee}}(t,x)$ which satisfies 
equation \eqref{7} taking $v^{\ee}$ instead of $v$.
Notice that equation \eqref{P0} has an additional term, in comparison with equation \eqref{4} which is given by the path-wise integral
\begin{equation*}
\int_0^{t}\,\langle\mathbf{G}_{\ubar{\al},\ubar{\del}}(t-s,x-\cdot)\sigma(u^{\ee,v^{\ee}}(s,
\cdot)),v^{\ee}(s,\cdot)\rangle_{\cH}\,ds,
\end{equation*}
however, the estimates \eqref{P1}, \eqref{P2}, \eqref{P-3} holds true and we proceed as follow.\medskip\\
\noindent As in \cite[Remark 2.2]{S-L}, $L^q(\Om)$ estimates of the 
the first and second terms in the right hand side of equation \eqref{7} leads,
up to a constant, to the same upper bound. Indeed, since $\|v^{\ee}\|_{\Ht}\leq N$ a.s., 
Cauchy-Schwartz's inequality on the Hilbert space $\Ht$ yields, for $q\geq1$
 \begin{align*}
\E \left| \int_0^{t}\,\langle
\mathbf{G}_{\ubar{\al},\ubar{\del}}(t-s,x-\cdot)\sigma(u^{\ee,v^{\ee}}(s,\cdot)),v^{\ee}(s,\cdot)\rangle_{\cH}\,ds \right|^{q}\hspace*{6.5cm}\\
\hspace*{6.5cm}\leq N^q\cdot \E\left|\int_0^{t}\,\left\|
\mathbf{G}_{\ubar{\al},\ubar{\del}}(t-s,x-\cdot)\sigma(u^{\ee,v^{\ee}}(s,\cdot))\right\|_{\cH}^2\,ds\right|^{\frac{q}{2}},
\end{align*}
Now, by using Burkholder's inequality to the stochastic integral we obtain 
\begin{align*}
    \E\left|\int_0^{t}\,\langle\mathbf{G}_{\ubar{\al},\ubar{\del}}(t-s,x-\cdot)\sigma(u^{\ee,v^{\ee}}(s,
\cdot)),e_k\rangle_{\cH}\,dB_s^k\right|^q\hspace*{7.05cm}\\
\hspace*{7cm}\leq C\cdot\E\left|\int_0^{t}\,\left\|\mathbf{G}_{\ubar{\al},\ubar{\del}}(t-s,x-\cdot)\sigma(u^{\ee,v^{\ee}}(s,
\cdot))\right\|_{\cH}^2\,ds\right|^\frac{q}{2},
\end{align*}
and this yields the extension to Theorem 2.1 in \cite{B-E-M} to cover equation \eqref{7}.
\end{dem}
\begin{rmk}\label{RM}
 The question of existence and uniqueness to the deterministic evolution equation 
 defined by \eqref{6} will be a straightforward consequence of the last Proposition taking $\ee=0$ in \eqref{7}.
\end{rmk}
\newpage%%
\noindent The next proposition is devoted to check the H\"{o}lder regularity of the stochastic
integral with respect to the martingale measure $F$. For the proof we refer the reader 
to Proposition 3.2 in \cite{B-E-M}.
\begin{prop}
\label{13}
 Let $\{U_{\ee,v^{\ee}}(t,x);\,(t,x)\in[0,T]\times\red\}$ be the stochastic integral
with respect to the martingale measure $F$ given by
$$U_{\ee,v^{\ee}}(t,x)=\sum_{k\geq1}\int_0^t\,\langle\mathbf{G}_{\ubar{\al},\ubar{\del}}(t-s,x-\cdot)\sigma(u^{\ee,v^{\ee}
}(s,\cdot)),e_k\rangle_{\cH}\,dB_s^k.$$
Then, under $\mathbf{(C)}$ and $\mathbf{(H_{\eta}^{\ubar{\al} })}$, $\eta\in(0,1]$, we have
\begin{enumerate}
\item[i)] For each $x\in\red$, a.s. $t\longmapsto U_{\ee,v^{\ee}}(t,x)$ is
$\beta_1$-H\"{o}lder continuous for $\beta_1\in(0,\frac{\al_0(1-\eta)}{2}),$
\item[ii)]For each $t\in[0,T]$, a.s. $x\longmapsto U_{\ee,v^{\ee}}(t,x)$ is
$\beta_2$-H\"{o}lder continuous for $\beta_2\in(0,\min\{1-\eta,\frac{1}{2}\}).$ 
\end{enumerate}
where $\al_0=\min_{1\leq i\leq d}\{\al_i\}$.
\end{prop}
Now, we give the H\"older regularity to the controlled equation \eqref{7}.
\begin{prop}\label{25}
Assume that $\mathbf{(C)}$ and $\mathbf{(H_{\eta}^{\ubar{\al }})}$, $\eta\in(0,1]$, holds and let
$u^{\ee,v^{\ee}}$ be the solution to equation \eqref{7}. Then 
$u^{\ee,v^{\ee}}$ belongs a.s. to the space ${\cal E}^{\beta}$ of $(\beta_1,\beta_2)$--H\"older continuous functions 
in time and space respectively. That is, for $(t',x')\neq(t,x)\in[0,T]\times K$
\begin{equation}\label{reg}
 \E\left(\left| u^{\ee,v^{\ee}}(t',x')- u^{\ee,v^{\ee}}(t,x)\right|^q\right)\leq
 C_q\left[\left|t'-t\right|^{q\cdot\beta_1}+\left\|x'-x\right\|^{q\cdot\beta_2}\right],
\end{equation}
$K$ being a compact subset of $\red$.
Moreover, for any $q\in[2,\infty[$
\begin{equation}\label{26}
 \sup_{\ee\leq1}\sup_{v^{\ee}\in\mathcal{A}_{\cH}^N}\E\left\|u^{\ee,v^{\ee}} \right\|_{\beta,K}^q<\infty.
\end{equation}
\end{prop}
\begin{dem}
 For any $(t',x'),(t,x)\in[0,T]\times K$ such that $t'\neq t$ and $x'\neq x$,
and for $q\in[2,\infty[$, consider $u^{\ee,v^{\ee}}$ and $Z^v$ the solution 
to equations \eqref{7} and \eqref{6} respectively,  then we have
\begin{align*}
\E\left(\left|\left[u^{\ee,v^{\ee}}(t',x')-Z^{v}(t',x')\right]-\left[u^{\ee,v^{\ee}}(t,x)-Z^{v}(t,x)\right]\right|^q\right)
&\leq C_q \E\left(\left|u^{\ee,v^{\ee}}(t',x')-u^{\ee,v^{\ee}}(t,x)\right|^q\right)\\
&+ C_q \E\left(\left|Z^{v}(t',x')-Z^{v}(t,x)\right|^q\right).\\
\end{align*}
Thus, the estimates on the increments \eqref{10} will be a consequence of \eqref{reg}
since $Z^v$ is a particular case of $u^{\ee,v^{\ee}}$ taking $\ee=0$ in equation \eqref{7}.
\noindent Now, let us focus on proving \eqref{reg}.
\begin{align*}
 \E\left(\left| u^{\ee,v^{\ee}}(t',x')- u^{\ee,v^{\ee}}(t,x)\right|^q\right)
&\leq2^{2q-2}\E\left|\sqrt{\ee}\sum_{k\geq1}\int_0^{t'}\,\langle\mathbf{G}_{\ubar{\al},\ubar{\del}}(t'-s,
x'-\cdot)\sigma(u^{\ee,v^{\ee}}(s,\cdot)),e_k\rangle_{\cH}\,dB_s^k\right.\\
&-\left.\sqrt{\ee}\sum_{k\geq1}\int_0^{t}\,\langle\mathbf{G}_{\ubar{\al},\ubar{\del}}(t-s,x-\cdot)\sigma(u^{
\ee,v^{\ee}}(s,\cdot)),e_k\rangle_{\cH}\,dB_s^k\right|^q\\
&+2^{2q-2}\E\left|\int_0^{t'}\,\langle\mathbf{G}_{\ubar{\al},\ubar{\del}}(t'-s,x'-\cdot)\sigma(u^{\ee,v^{\ee}}(s,
\cdot)),v^{\ee}(s,\cdot)\rangle_{\cH}\,ds\right.\\
&-\left.\int_0^{t}\,\langle\mathbf{G}_{\ubar{\al},\ubar{\del}}(t-s,x-\cdot)\sigma(u^{\ee,v^{\ee}}(s,\cdot)),v^{\ee}(s,
\cdot)\rangle_{\cH}\,ds\right|^q\\
&+2^{2q-2}\E\left|\int_0^{t'}\int_{\red}\mathbf{G}_{\ubar{\al},\ubar{\del}}(t'-s,x'-y)
b(u^{\ee,v^{\ee}}(s,y))ds\, dy\right.\\
&-\left.\int_0^t\int_{\red}\mathbf{G}_{\ubar{\al},\ubar{\del}}(t-s,x-y) b(u^{\ee,v^{\ee}}(s,y))ds\, dy\right|^q\\
&=C_q\sum_{i=1}^3 \Lambda_i,
\end{align*}
where
 \begin{align*}
\Lambda_1&=\E\left|\sqrt{\ee}\sum_{k\geq1}\int_0^{t'}\,\langle\mathbf{G}_{\ubar{\al},\ubar{\del}}(t'-s,x'-\cdot)\sigma(u^{\ee,v^{\ee}}(s,\cdot)),e_k\rangle_{\cH}\,dB_s^k\right.\\
&-\left.\sqrt{\ee}\sum_{k\geq1}\int_0^{t}\,\langle\mathbf{G}_{\ubar{\al},\ubar{\del}}(t-s,x-\cdot)\sigma(u^{\ee,v^{\ee}}(s,\cdot)),e_k\rangle_{\cH}\,dB_s^k\right|^q,\\
\Lambda_2&=\E\left|\int_0^{t'}\,\langle\mathbf{G}_{\ubar{\al},\ubar{\del}}(t'-s,x'-\cdot)\sigma(u^{\ee,v^{\ee}}(s,\cdot)),v^{\ee}(s,\cdot)\rangle_{\cH}\,ds\right.\\
&-\left.\int_0^{t}\,\langle\mathbf{G}_{\ubar{\al},\ubar{\del}}(t-s,x-\cdot)\sigma(u^{\ee,v^{\ee}}(s,\cdot)),v^{\ee}(s,
\cdot)\rangle_{\cH}\,ds\right|^q,
\intertext{and}
\Lambda_3&=\E\left|\int_0^{t'}\int_{\red}\mathbf{G}_{\ubar{\al},\ubar{\del}}(t'-s,x'-y)
b(u^{\ee,v^{\ee}}(s,y))ds\, dy\right.\\
&-\left.\int_0^t\int_{\red}\mathbf{G}_{\ubar{\al},\ubar{\del}}(t-s,x-y) b(u^{\ee,v^{\ee}}(s,y))ds\, dy\right|^q.
\end{align*}
\noindent As mentioned before in the proof of Proposition~\ref{12}, up to a constant, $\Lambda_1$ and
$\Lambda_2$ have the same upper bound, which is, from Proposition~\ref{13}, given by
 $$C_q\left[\left|t'-t\right|^{q\cdot\beta_1}+\left\|x'-x\right\|^{q\cdot\beta_2}\right].$$
Now, after a change of variable, $\Lambda_3$ becomes
\begin{align*}
\Lambda_3=\E&\left|\int_0^{t}\int_{\red}\mathbf{G}_{\ubar{\al},\ubar{\del}}(t-s,x-y)\right.\times\left[b\left(u^{\ee,v^{\ee}}(s+t'-t,y+x'-x)\right)
-b\left(u^{\ee,v^{\ee}}(s,y)\right)\right]ds\, dy\\
&-\left.\int_0^{t'-t}\int_{\red}\mathbf{G}_{\ubar{\al},\ubar{\del}}(t'-s,x'-y) b(u^{\ee,v^{\ee}}(s,y))ds\, dy\right|^q,
\end{align*}
and we proceed as in the proof of \cite[Theorem 3.1]{B-E-M} setting $h=t'-t$ and
$z=x'-x$. That is, H\"older's inequality, assertion $(i)$ of Lemma \ref{lqr} along with the Lipschitz condition and linear growth property of $b$ imply
\begin{align*}
 \Lambda_3&\leq C_q\int_0^{t}\int_{\red}\mathbf{G}_{\ubar{\al},\ubar{\del}}(t-s,x-y)\times\E\left|b\left(u^{\ee,v^{\ee}}(s+t'-t,y+x'-x)\right)
 - b\left(u^{\ee,v^{\ee}}(s,y)\right)\right|^q dsdy\\
&\qquad\qquad\qquad+\int_0^{t'-t}\int_{\red}\mathbf{G}_{\ubar{\al},\ubar{\del}}(t'-s,x'-y)\E\left|b(u^{\ee,v^{\ee}}(s,y))\right|^q dsdy\\
&\leq C_q\left[|t'-t|+\int_0^{t}\sup_{y\in\red}\E\left|u^{\ee,v^{\ee}}
(s+t'-t,y+x'-x)- u^{\ee,v^{\ee}}(s,y)\right|^qds\right].
\end{align*}
Putting together all the estimates and using Gronwall's lemma, we conclude the
proof by the Kolmogorov continuity criterium.\medskip\\
Notice that going through the arguments, we can easily get uniform estimates for $u^{\ee,v^{\ee}}$ in $\ee\in]0,1]$ and $v^{\ee}\in\mathcal{A}_{\cH}^N$,
therefore \eqref{26} remain valid.
\end{dem}
\begin{prop}\label{27}
 Assuming $\mathbf{(C)}$ and $\mathbf{(H_{\eta}^{\ubar{\al }})}$, $\eta\in(0,1]$, let
$\{v,v^{\ee};\,\ee>0\}\subset\mathcal{A}_{\cH}^N,$ such that
${P}.\,a.s.$
\beq\label{00}
\lim_{\ee\rightarrow0}\,\langle v^{\ee}-v,\,g\,\rangle_{\Ht}=0,\,\,
\textit{for any }\,g\in\Ht.
\eeq
Then, for any $(t,x)\in[0,T]\times K$, $q\in[2,\infty[$ we have
$$\lim_{\ee\rightarrow0}\E\left(\left|u^{\ee,v^{\ee}}(t,x)-Z^v(t,
x)\right|^q\right)=0.$$
\end{prop}
\begin{dem}
First we need to recall the following two key ingredients
\begin{equation}\label{20}
 \int_0^t ds\int_{\red}\left|\mathcal{F}\mathbf{G}_{\ubar{\al},\ubar{\del}}(t-s,x-\cdot)(\xi)\right|^2\mu(d\xi)<\infty,
\end{equation}
and
\beq\label{30}
\sup_{v\in\mathcal{A}_{\cH}^N}\sup_{(t,x)\in[0,T]\times\red}
\E\left[\left|Z^v(t,x)\right|^q\right]<\infty.
\eeq
\noindent Since condition $\mathbf{(H_{\eta}^{\ubar{\al }})}$ holds for $\eta\in(0,1]$, then \eqref{20} follow.
Now, from the fact that $Z^v$ is the solution to the particular equation \eqref{7} taking $\ee=0$,
then \eqref{30} is an immediate consequence of \eqref{estm}.
\noindent Fix $q\in[2,\infty[$, then
$$\E\left(\left|u^{\ee,v}(t,x)-Z^v(t,x)\right|^q\right)\leq C_q\sum_{i=1}^4\E\left|A_{i,\ee}(t,x)\right|^q,$$ 
where 
\begin{align*}
&A_{1,\ee}(t,x)=\int_0^{t}\int_{\red}\mathbf{G}_{\ubar{\al},\ubar{\del}}(t-s,x-y)\left[b(u^{\ee,v^{\ee}}(s,y))-b(Z^{v}(s,y))\right] ds\, dy,\\ 
&A_{2,\ee}(t,x)=\sqrt{\ee}\sum_{k\geq1}\int_0^{t}\,\langle\mathbf{G}_{\ubar{\al},\ubar{\del}}(t-s,x-\cdot)\sigma(u^{\ee,v^{\ee}}(s,\cdot)),e_k\rangle_{\cH}\,dB_s^k,\\
&A_{3,\ee}(t,x)=\int_0^{t}\,\langle\mathbf{G}_{\ubar{\al},\ubar{\del}}(t-s,x-\cdot)\left[\sigma(u^{\ee,v^{\ee}}(s,\cdot))-
\sigma(Z^v(s,\cdot))\right],v^{\ee}(s,\cdot)\rangle_{\cH}\,ds,\\
&A_{4,\ee}(t,x)=\int_0^{t}\,\langle\mathbf{G}_{\ubar{\al},\ubar{\del}}(t-s,x-\cdot)
\sigma(Z^v(s,\cdot)),v^{\ee}(s,\cdot)-v(s,\cdot)\rangle_{\cH}\,ds.
\end{align*}
For the first term $A_{1,\ee}$, by H\"older's inequality along with the Lipschitz condition on $b$ we get
\begin{align*}
\E\left|A_{1,\ee}(t,x)\right|^q &\leq 
\int_0^{t}\int_{\red}\mathbf{G}_{\ubar{\al},\ubar{\del}}(t-s,x-\cdot)\E\left|b(u^{\ee,v^{\ee}}(s,y))-b(Z^{v}(s,y))\right|^q ds\, dy\\
&\leq C_q\int_0^{t}\sup_{(r,y)\in[0,s]\times\red}\E\left|u^{\ee,v^{\ee}}(r,y)-Z^v(r,y)\right|^q ds.
\end{align*}
For the second term $A_{2,\ee}$, Burkholder's inequality together with the linear growth property of $\sigma$ and \eqref{30}
 yields
\begin{align*}
\E\left|A_{2,\ee}(t,x)\right|^q &=\ee^{\frac{q}{2}}\E\left(\int_0^{t}\,\left\|\mathbf{G}_{\ubar{\al},\ubar{\del}}(t-s,x-\cdot)\sigma(u^{\ee,v^{\ee}}(s,
\cdot))\right\|_{\cH}^2\,ds\right)^\frac{q}{2}\\
&\leq\ee^{\frac{q}{2}}\int_0^tds\left(1+\sup_{(r,y)\in[0,s]\times\red}\E\left|u^{\ee,v^{\ee}}(r,y)\right|^q\right)\times(\mathcal{J}(t-s))\\
&\times  \left(\int_0^t ds\int_{\red}\left|\mathcal{F}\mathbf{G}_{\ubar{\al},\ubar{\del}}(t-s,x-\cdot)(\xi)\right|^2\mu(d\xi)\right)^{\frac{q}{2}-1}\\
&\leq C_q\ee^{\frac{q}{2}},
\end{align*}
with ${\cal J}(t-s)$ given by \eqref{J}.\medskip\\
Then, we have $$\lim_{\ee\rightarrow0}\E\left|A_{2,\ee}(t,x)\right|^q=0.$$
To deal with the term $A_{3,\ee}$, first we apply the Cauchy-Schwartz's inequality to the
 inner product on $\cH$, and the property $\sup_{\ee\leq1}\|v^{\ee}\|_{\Ht}\leq N$ we have
\begin{align*}
\E\left|A_{3,\ee}(t,x)\right|^q &\leq\E\left(\int_0^{t}\left\|\mathbf{G}_{\ubar{\al},\ubar{\del}}(t-s,x-\cdot)\left[\sigma(u^{\ee,v^{\ee}}(s,\cdot))-
\sigma(Z^v(s,\cdot))\right]\right\|^2_{\cH}ds\right)^{\frac{q}{2}}\\
&\times\left(\int_0^{t}\left\|v^{\ee}(s,\cdot)\right\|^2_{\cH}\right)^{\frac{q}{2}}\\
&\leq C_q\E\left(\int_0^{t}\left\|\mathbf{G}_{\ubar{\al},\ubar{\del}}(t-s,x-\cdot)\left[\sigma(u^{\ee,v^{\ee}}(s,\cdot))-
\sigma(Z^v(s,\cdot))\right]\right\|^2_{\cH}ds\right)^{\frac{q}{2}}.
\end{align*}
Next, H\"older's inequality with respect to the measure on $[0,T]\times\red$ given by $|\mathcal{F}\mathbf{G}_{\ubar{\al},\ubar{\del}}(t-s)(\xi)|^2\mu(d\xi)ds$,
 \eqref{20} and the Lipschitz condition on $\sigma$ yields
\begin{align*}
 \E&\left|A_{3,\ee}(t,x)\right|^q \leq C_q\left(\int_0^t ds\int_{\red}\left|\mathcal{F}\mathbf{G}_{\ubar{\al},\ubar{\del}}(t-s,x-\cdot)
(\xi)\right|^2\mu(d\xi)\right)^{\frac{q}{2}-1} \\
&\times\int_0^{t}\sup_{(r,y)\in[0,s]\times\red}\E\left|u^{\ee,v^{\ee}}(r,y)-Z^v(r,y)\right|^q
\left({\cal J}(t-s)\right)ds\\
&\leq C_q \int_0^{t}\sup_{(r,y)\in[0,s]\times\red}\E\left|u^{\ee,v^{\ee}}(r,y)-Z^v(r,y)\right|^q\left({\cal J}(t-s)\right) ds,
\end{align*}
where ${\cal J}(t-s)$ is defined by \eqref{J}.\medskip\\
\noindent For the last term, we apply Cauchy-Schwartz's inequality to the inner product on $\cH$, H\"older's inequality with respect to the measure 
 $|\mathcal{F}\mathbf{G}_{\ubar{\al},\ubar{\del}}(t-s)(\xi)|^2\mu(d\xi)ds$ on $[0,T]\times\red$, the linear growth property of $\sigma$ 
 and \eqref{30}, then we get
\begin{align*}
 \E\left|A_{4,\ee}(t,x)\right|^q &\leq \E\left(\int_0^t\left\|\mathbf{G}_{\ubar{\al},\ubar{\del}}(t-s,x-\cdot)
 \sigma(Z^v(s,\cdot))\right\|^2_{\cH}ds\right)^{\frac{q}{2}}\\
 &\times\left(\int_0^t\left\|v^{\ee}(s,\cdot)-v(s,\cdot)\right\|_{\cH}^2ds\right)^{\frac{q}{2}}\\
 &\leq C_q\int_0^tds\left(1+\sup_{(r,y)\in[0,s]\times\red}\E\left|Z^v(r,y)\right|^q\right)\times(\mathcal{J}(t-s))\\
&\times  \left(\int_0^t ds\int_{\red}\left|\mathcal{F}\mathbf{G}_{\ubar{\al},\ubar{\del}}(t-s,x-\cdot)(\xi)\right|^2\mu(d\xi)\right)^{\frac{q}{2}-1}\\
&\times\left(\int_0^t\left\|v^{\ee}(s,\cdot)-v(s,\cdot)\right\|_{\cH}^2ds\right)^{\frac{q}{2}}\\
 &\leq C_q\left\|v^{\ee}-v\right\|_{\Ht}^{q}.
\end{align*}
Thus, as $\ee$ goes to $0$, \eqref{00} imply
$$\lim_{\ee\rightarrow0}\E\left|A_{4,\ee}(t,x)\right|^q=0.$$
Now, let
$$\Phi_{\ee}(t)=\sup_{(t,x)\in[0,T]\times\red}\E\left(\left|u^{\ee,v^{\ee}}(t,x)-Z^v(t,x)\right|^q\right).$$
Then, taking together all the estimates, we get
$$\Phi_{\ee}(t)\leq C_q\left[\ee^{\frac{q}{2}}+\E\left|A_{4,\ee}(t,x)\right|^q+\int_0^t\Phi_{\ee}(s)\left(1+{\cal J}(t-s)\right)ds\right],$$
${\cal J}(t-s)$ defined by \eqref{J}.
\newpage
 We end the proof by applying the extended version of Gronwall's Lemma in \cite[Lemma 15]{Da}, and we find
\begin{equation*}
 \lim_{\ee\rightarrow0}\sup_{(t,x)\in[0,T]\times\red}\E\left(\left|u^{\ee,v^{\ee}}(t,x)-Z^v(t,x)\right|^q\right)=0.
\end{equation*}
\end{dem}
\newline
\begin{dem}(of Theorem \ref{21})
\\%
 By Proposition \ref{25} and Proposition \ref{27}, the estimation on increments \eqref{10} 
and Point-wise convergence \eqref{11} holds true. 
Then, as it has been argued before, Theorem \ref{21} will follow.
\end{dem}
%%%%%%%%%%%%%%%%%%%%%%%%%%%%%%%%%%%%%%%%%%%%%%%%%%%%%%%%%%%  BIBLIOGRAPHY  %%%%%%%%%%%%%%%%%%%%%%%%%%%%%%%%%%%%%%%%%%%%%%%%%%%%%%%%%%%

\vspace*{1.5cm}
{\scshape
\begin{flushright}
\begin{tabular}{l}
Universit\'{e} Cadi Ayyad\\
Facult\'e des Sciences Semlalia\\
D\'{e}partement des Math\'{e}matiques\\
LIBMA, B.P. 2390, Marrakech \\
Maroc \\
{\upshape e-mail: \href{mailto: t.elmellali@ced.uca.ma}{\nolinkurl{t.elmellali@ced.uca.ma}},}\\
\\
Universit\'{e} Paris Descartes \\
MAP5, CNRS UMR 8145 \\
45, rue des Saints-P\`eres \\
75270 Paris Cedex 6 \\
France \\
{\upshape e-mail: \href{mailto: mohamed.mellouk@parisdescartes.fr}{\nolinkurl{mohamed.mellouk@parisdescartes.fr}}}\\
\end{tabular}
\end{flushright}
}
\end{document}